\documentclass[10pt]{amsart}
\usepackage{caption}
\usepackage{subcaption}
\usepackage{enumerate,amssymb,  mathrsfs, graphicx}
\newtheorem{theorem}{Theorem}[section]

\newtheorem*{lemma*}{Lemma}

\theoremstyle{definition}

\newtheorem{remark}[theorem]{Remark}

\theoremstyle{remark}

\numberwithin{equation}{section}

\newcommand{\A}{\mathbb{A}}
\newcommand{\B}{\mathbb{B}}

\newcommand{\E}{\mathcal{E}}

\newcommand{\X}{\mathbb{X}}

\newcommand{\Y}{\mathbb{Y}}

\begin{document}
\title[Minimization of Dirichlet energy of $j-$degree mappings]{Minimization of Dirichlet energy of $j-$degree mappings between annuli}

\author[D. Kalaj]{David Kalaj}
\address{University of Montenegro, Faculty of natural sciences and mathematics, Podgorica, Cetinjski put b.b. 81000 Podgorica, Montenegro }
\email{davidk@ucg.ac.me}

%    General info
\subjclass[2010]{Primary 35J60; Secondary  30C70}

%\date{\today}

\keywords{$j-$degree mappings, harmonic mappings, Dirichlet energy. }

%\begin{abstract}
%\end{abstract}

\maketitle

%\dedicatory{}
\begin{abstract}
Let $\mathbb{A}$ and $\B$ be circular annuli in the complex plane and consider the Dirichlet energy integral of $j-$degree mappings between $\A$ and $\B$. Then we minimize this energy integral. The minimizer is a $j-$degree harmonic mapping between annuli $\A$ and $\B$ provided it exits. If such a harmonic mapping does not exist, then the minimizer is still a $j-$degree mapping which is harmonic in $\A'\subset \A$ and it is a squeezing mapping in its complementary annulus $\A''=\A\setminus \A$. Such a result is an extension of the certain result of Astala, Iwaniec and Martin \cite{astala2010}.
\end{abstract}

%\maketitle \pagestyle{myheadings} \markboth{ David Kalaj}{Dirichlet-type energy of mappings between two concentric annuli}
\section{Introduction}\label{sec-1}
In this paper, we continue to study the minimization problem of the {\it Dirichlet energy} of Sobolev mappings belonging to the  $W^{1,2}$ class,  between domains in $\mathbb{R}^2$. Such research is related to the principle of non-interpenetration of matter in the mathematical theory
of Nonlinear Elasticity (NE) see for example \cite{sverak}. We will minimize the energy of $j-$degree mappings between certain doubly connected domains.

\subsection{Dirichlet integral} Let $D$ and $\Omega$ be two bounded planar domains in $\mathbb{R}^2$, and let $h:D\stackrel{\text{onto}}{\longrightarrow}\Omega$ be a mapping that belongs to the Sobolev class $W^{1,2}$.
The {\it Dirichlet integral} (also called {\it conformal energy}) of $h$ is given as follows:
$$\mathcal{E}[h]=\int_{D}|Dh(z)|^2dxdy,$$
where $z=(x, y)\in D$ and $|Dh(z)|^2 =|f_z|^2+|f_{\bar z}|^2=\frac{1}{2}(|f_x|^2+|f_y|^2)$.

Let $j$ be a positive integer and assume that $\mathcal{H}_j(D,\Omega)$ is the class of $j-$degree smooth mappings  between $D$ and $\Omega$ mapping the inner/outer boundary onto inner/outer boundary. Recall that we say that a mapping $f$ has degree $j$ in a regular point $y\in \Omega$ if $$j=\deg(f,\Omega,y)=\sum_{x\in f^{-1}(y)}\mathrm{sign}(\det (Df(x))).$$ Here regular means that $\det (Df(x))\neq 0$ for $x\in f^{-1}(y)$. For non-regular points, the degree is defined throughout a sequence of regular points $y_n$ converging to $y$.  This sequence exists according to the Sard theorem. Then we say that the mapping has a $j-$degree if $g$ has degree $j$ in every point from the image domain. The notation of degree can be also extended to Sobolev mappings of $W^{1,2}$ class, or even to the continuous mappings.

If $D$ and $\Omega$ are two simply connected domains in $\mathbb{R}^2$ different from $\mathbb{R}^2$ and let $\mathbb{D}$ be the unit disk. Then the Riemann mapping theorem asserts that there exist conformal bijections  $g:D \stackrel{\text{onto}}{\longrightarrow}  \mathbb{D}$ and $k: \mathbb{D}\stackrel{\text{onto}}{\longrightarrow} \Omega$. Then the mapping $f(z) = k(g^j(z))$ is $j-$defree conformal mapping between $D$ and $\Omega$.

Moreover, this map is a minimizer of the Dirichlet integral throughout the $j-$degree mappings. Indeed, for every $h$ we have
$$\mathcal{E}[h]=\int_{D}|Dh(z)|^2dxdy\geq \int_{D} \det(Dh(z))dxdy=j\mbox{Area}(\Omega).$$
The equality holds if and only if $h$ is $j-$degree  conformal (because in this case, one has $|Dh|^2=\frac{1}{2}\left(|h_x|^2+|h_y|^2\right)=\det(Dh(z))$).

If $D$ and $\Omega$ are doubly connected domains, then the situation is very different.
Namely, if the doubly connected domains are not conformally equivalent, then such a conformal minimizer does not exist. However, in this case, the harmonic diffeomorphisms can be an essential replacement for the conformal diffeomorphisms. The existence problem has been studied in a large number of papers. See for example the references  \cite{astala2010, inventiones, calculus, klondon,  koh1}, which deal with the existence of minimisers  for two-dimensional domains and surfaces and with their properties. The authors of \cite{li} considered a similar problem for the so-called $\sigma_2$ energy between annuli in $\mathbb{R}^4$.

The problems of this kind are related to the so called {\it Nitsche conjecture} for harmonic mappings.
Let us briefly explain what is its formulation.

\subsection{Nitsche conjecture} Let $\A:=\{z:r<|z|<R\}$ and $\B=\{w:r_\ast<|w|<R_\ast\}$ be two annuli.
In 1962 J.C.C. Nitsche announced, in a short article ~\cite{NCONJ}, that the existence of a harmonic homeomorphism
$h \colon \A  \overset{\textnormal{\tiny{onto}}}{\longrightarrow} \B$, whether or not it comes from a minimal graph, yields a lower bound on $\mathrm{Mod}\,(\B)$ in terms of $\mathrm{Mod}\,(\A)$. He conjectured that the necessary and sufficient condition for  such a mapping to exist is  the following inequality, now known as the {\it Nitsche bound}
\begin{equation}\label{nb}
\frac{R_\ast}{r_\ast} \ge \frac{1}{2} \left(\frac{R}{r}+ \frac{r}{R}\right).
\end{equation}
 Various lower bounds for $R_\ast/r_\ast$ have been obtained by Lyzzaik~\cite{Lyz}, Weitsman~\cite{weit}, Kalaj~\cite{Ka}.
  Finally, this conjecture was solved by Kovalev, Onninen, and Iwaniec \cite{solut} who give an affirmative answer to it.

\subsection{Motivations}
An indirect evidence of Nitsche conjecture has been given previously by Astala, Iwaniec and Martin in \cite{astala2010}. Namely they showed that the minimizer of Dirichlet energy of $W^{1, 2}$- of homeomorphisms exists and present itself a $W^{1, 2}$- of homeomorphism precisely when  the Nitsche bound \eqref{nb} is satisfied.

On the other hand side, assume that reference annulus $\A$ is substantially fatter than $\A_*$, precisely
\begin{equation}\label{21-10-9}
\frac{R}{r}>\frac{R_\ast}{r_\ast}+\sqrt{\frac{R^2_\ast}{r^2_\ast}-1}, \ \ \ \mbox{or equivalently}\ \ \ \frac{R_\ast}{r_\ast}<\frac{1}{2}\left(\frac{R}{r}+\frac{r}{R}\right).
\end{equation}
Let $\rho$ be determined by the so-called critical Nitsche equation
$$\frac{R_\ast}{r_\ast}=\frac{1}{2}\left(\frac{R}{\rho}+\frac{\rho}{R}\right).$$
Then in \cite{astala2010} it is showed that the minimizer within all weak $W^{1, 2}$-limits of homeomorphisms: $h:\A\stackrel{\text{onto}}{\longrightarrow}\A^\ast$
takes the form:
$$h(z)=\left\{
  \begin{array}{ll}
  \vspace{2mm}
    r_\ast\frac{z}{|z|}, &r<|z|<\rho\ \ \   \hbox{-not harmonic, squeezing map} \\
    r_\ast\left(\frac{z}{2\rho}+\frac{\rho}{2\bar{z}}\right), & \rho<|z|<R\ \ \ \hbox{-critical harmonic Nitsche map}.
  \end{array}
\right.$$
Note that $\det(Dh(z))=|h_z|^2-|h_{\bar z}|^2= \frac{r_\ast^2}{4|z|^2}-\frac{r_\ast^2}{4|z|^2}\equiv 0$ for $r\leq|z|\leq \rho$. This minimizer is unique up to the rotation of annuli as it is shown in \cite{astala2010}.
\subsection{Statement of the problem and the formulation of the main result}
Let $\X$ and $\Y$ be two doubly connected domains in $\mathbb{R}^2$.

This paper aims  to minimize the Dirichlet energy integral $\mathcal{E}[h]$ for mappings $h$ belonging to the following class:
\begin{align*}
  \mathcal{H}_j(\X,\Y):= & \big\{\mbox{the $j-$degree mappings in the Sobolev class}\  {W}^{1,2}(\X, \Y), \\
   & \mbox{which map the inner/outer boundary onto the inner/outer boundary}\big\}.
\end{align*}
In this paper we assume that $\Y=\B(r^\ast, R^\ast)$.
It should be noted that for every doubly connected domain $\X$  there is a conformal diffeomorphism $\zeta: \X \mapsto \A(r,R)$. In this case, the conformal modulus of $\X$ is defined by $\mathrm{Mod}(\X)=\log \frac{R}{r}$.  Without losing of generality we assume that $\A=\A(1,r)$ and  $\B=A(1,R)$, where $r,R>1$, because the minimizations of Dirichlet energy of  mappings $g$ and their transformations $h(z) = \alpha g(\beta z)$ is equivalent.

The following theorem is the main result of this paper.

\begin{theorem}\label{theoj}
There exists a $j-$degree radial harmonic mapping between $\A$ and $\B$ if and only if the condition \eqref{Nitchegeneral} is satisfied. In this case, the harmonic mapping $g^\circ: \A\to \B$  is given by \eqref{cini}.

Moreover the minimum of the Dirichlet energy \begin{equation}\label{diri}\E[f]=\int_{\A} |Df|^2 dxdy,\end{equation} for $f\in \mathcal{H}_j(\A,\B)$ is attained for the $j-$degree radial harmonic mapping $g^\circ$ provided that \eqref{Nitchegeneral} is satisfied. If \eqref{Nitchegeneral} is not satisfied, then the minimum of \eqref{diri} is attained for a $j-$degree mapping  $g^\diamond$ (see \eqref{hcirc} below) which is harmonic in a sub-annulus of $\A$ but it is an "squeezing" in rest of $\A$.
The minimum on both cases is unique up to a rotation.

\end{theorem}

\begin{remark}
Theorem~\ref{theoj} is a generalization of the similar result obtained by Astala, Iwaniec, and Martin in \cite{astala2010} for smooth homeomorphismsm.  It must be emphasized that Theorem~\ref{theoj} however cannot be deduced from such a result, because $j-$degree mappings cannot be decomposed as the composition of a homeomorphism and the power function $z^j$ as one can guess. It seems that such a result is new even for degree one mappings, which are not necessarily homeomorphisms.   The proof is an adaptations of the corresponding proofs in \cite{ka2019, arma2009} by employing the free Lagrangian for $j-$degree mappings mapping the inner/outer boundary to inner/outer boundary.
\end{remark}

\section{Radial harmonic mappings and radial minimizers}

Set $f(z) = e^{j \imath \tau} G(\rho)$, and call this mapping \emph{$j-$degree radial mapping} where $z=t e^{\imath \tau}$ and solve the Laplace equation $\Delta f=0$. Then we get the second order differential equation \begin{equation}\label{firsteq}-{j^2 G(t)}+t{G'(t)}+t^2G''(t)=0.\end{equation}
By taking the substitution $t= e^x$, $H(x)=G(e^x)$, we obtain $H'(x) = e^x G'(e^x)=t G'(t)$ and $$H''(x) =H'(x) + e^{2x} G''(e^x)=F'(x)+t^2 G''(t).$$
Thus \eqref{firsteq} is reduced to \begin{equation}\label{bulla}-j^2 H(x) +H''(x)=0.\end{equation}
Now take
$s=H(x)$, $H'(x)=h(s)$. Then we obtain $H''(x)= h'(s) h(s)$. And \eqref{bulla} is reduced to the equation $$2j^2 s=(h^2(s))'.$$ Thus $$h(s) = \pm \sqrt{j^2 s^2+c}.$$
After returning back the variable $x=\log t$ and the function $G$,  we obtain the equation
\begin{equation}\label{pm}\pm  G'(t)\frac{1}{\sqrt{j^2G(t)^2+c_1}} +\frac{1}{t}=0,\end{equation} where plus sign is in the case when $G$ is decreasing and minus sign is for the case when $G$ is increasing.

 Then the general solution of \eqref{firsteq} is given by $$G(t)=a \cosh (j^2 \log(t))+b \sinh (j^2 \log (t)).$$

Assuming that it maps $\A(1,r)$ onto $\A(1,R)$, mapping the inner boundary onto the inner boundary and the outer boundary onto the outer boundary we get
\begin{equation}\label{cini}g^\circ(z) = \left(\frac{r^j R-r^{2 j}}{1-r^{2 j}}\right) \bar z^{-j}+\frac{\left(1-r^j R\right) z^j}{1-r^{2 j}}\end{equation}
and the given function \eqref{cini} maps $\A(1,r)$ onto $\A(1,R)$ if and only if it is satisfied the Nitsche type bound
\begin{equation}\label{Nitchegeneral}R\ge \frac{1}{2} r^{-j} \left(1+r^{2 j}\right).\end{equation}
Then the previous inequality can be written as $$r\le r_\circ :=\left(R+\sqrt{-1+R^2}\right)^{\frac{1}{j}}.$$
Indeed, $g^\circ(\rho)$ is increasing if and only if $$r^j \left(R-r^j\right)+\rho^{2 j} \left(r^j R-1\right)\ge 0,\ \ \ \rho\in[1,r],$$ and this inequality is equivalent with \eqref{Nitchegeneral}. Now if $r>r_\circ$, then there is an annulus $\A(\rho, r_\circ)$ having the conformal modulus equal to the modulus of  $\A(1,r)$. Indeed we have $\rho = r_\circ/r$.
Let \begin{equation}\label{hcirc}g^\diamond(z) = \left\{
                       \begin{array}{ll}
                         g^\circ(z), & \hbox{if $1\le |z|< r_\circ$;} \\
                         \frac{z^j}{|z|^j}, & \hbox{if $\rho< |z|\le 1$.}
                       \end{array}
                     \right.\end{equation}

This is a $j-$degree mapping. Namely for every $1<|w|<R$, there exist exact $j$ points from $1<|z|<r_\circ$ so that $g^\circ(z)=w$. Moreover the $\mathrm{sign}(\det(Dg^\circ(z)))=1$. We can also define the degree for points $|w|=1$ by using a sequence of points $1<|w_n|<R$ converging to $w$. See Figure~1 for a $2$-degree mapping $g^\diamond: \A(1/2,2)\to \A(1,17/8).$
\begin{figure}\label{figa}
\centering
\includegraphics{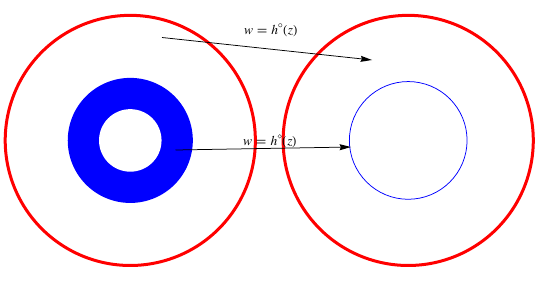}
\caption{The annulus $\A(1/2,2)=\A(1/2,1]\cup \A[1,2)$ is mapped onto the annulus $\A(1,17/8)$ in such a way that $\A(1/2,1)$ is squeezed in the unit circle and wrapped two times, while the other part is $2-$degree harmonic mapping of $\A(1,2)$ onto $\A(1,17/8)$. }
\label{21}
\end{figure}

\subsection{The energy of the stationary mapping}

Elementary calculations lead to
$$|Dg^\circ|^2=\frac{j^2 G^2}{t^2}+\dot G^2.$$
Thus it follows from (\ref{pm}) that
$$|Dg^\circ|^2=\frac{2j^2G^2+c_1}{t^2},$$ where $$c_1=\frac{4j^2 \left(1+R^2- R (r^j+r^{-j})\right)}{ (r^j-r^{-j})^2}.$$

Observe that $c_1=0$ if and only if $R=r^j$ or $ R=r^{-j}$. In this case $f(z) =z^{\pm j}$ is a conformal mapping between $\A(r)$ and $\A(R)$.

\begin{equation}\label{dgs}\dot G^2= \frac{(j^2)G^2+c_1}{t^2}, \end{equation}

Let
$F=G^{-1}$
then

\begin{equation}\label{fprims}\frac{F'(s)}{F(s)}=\frac{1}{\sqrt{j^2 s^2+c_1}}.\end{equation}

The energy of $g^\circ$ is now as follows
\begin{equation}\label{energyof}\mathcal{E}[g^\circ]=\int_{\A}|Dg^\circ|^2=4\pi j^2 \int_1^r\frac{G^2(t)}{t}dt+2\pi c_1\log r.\end{equation}
Now we calculate $\mathcal{E}[g^\diamond]$. In this case $c_1=-j^2$. We have \[\begin{split}\mathcal{E}[g^\diamond]&=\int_{\A(1,r)}|Dg^\circ|^2+\int_{\A(\rho ,1)}|Dg^\diamond|^2
\\&=4\pi j^2 \int_1^r\frac{G^2(t)}{t}dt+2\pi -j^2\log r+ 2\pi j^2\int_\rho^1 \frac{1}{t} dt\\&=4\pi j^2 \int_1^r\frac{G^2(t)}{t}dt-2\pi j^2\log \rho r.\end{split}\]

\section{Free Lagrangians}

Assume that $\mathcal{H}_j(\A,\B)$ is a class of $j-$degree mappings, mapping $\A$ onto $\B$ and mapping the inner/outer boundary to the inner/outer boundary.
Motivated by the paper of Iwaniec and Onninen \cite{iwon, arma2009, annalen2010} we consider the following free Lagrangians.

a)  A function in $t= |z|$;
\begin{equation}\label{functionx}L(z, h, Dh) dx\wedge dy = M(t) dx\wedge dy \end{equation}
Thus, for all $h \in  \mathcal{H}_j(\mathbb{A}, \mathbb{B})$ we have

\begin{equation}\label{simpli}\int_\A L(z, h, Dh) dx\wedge dy =
\int_\A M(|z|) dx\wedge dy=2\pi\int_r^R t M(t) dt.\end{equation}

b)  Pullback of a form in $\B$ via a given mapping $h \in  \mathcal{H}_j(\A, \B)$;
\begin{equation}\label{pullback}L(z, h, Dh) dx\wedge dy  = N(|h|) J (z, h) dx\wedge dy , \text{where} \ \ \ N \in  L^ 1(r,R) \end{equation}
Thus, for all $h \in  \mathcal{H}_j(\mathbb{A}, \mathbb{B})$ we have

\begin{equation}\label{easi0}\int_\A L(z, h, Dh) dx\wedge dy =j\int_{\B}
N(|w|) du\wedge dv = 2\pi j \int_{r_\ast}^{R_\ast}
sN(s )  ds . \end{equation}

c)  A radial free Lagrangian
\begin{equation}\label{radial}L(z, h, Dh) dx\wedge dy = A (|h|)\frac{|h|_N}{ |z|}
dx\wedge dy\end{equation} where $A \in  L^   1( r_\ast, R_\ast)$.
Thus, for all $h \in  \mathcal{H}_j(\mathbb{A}, \B)$ we have

\begin{equation}\label{easi}\int_\A L(z, h, Dh) dx\wedge dy = 2\pi \int_r^R A(|h|)\frac{\partial |h|}{\partial \rho} d\rho  = 2\pi\int_{r_\ast}^{R_\ast} A(s)ds  \end{equation}
d)  An angular free Lagrangian
\begin{equation}\label{angular}L(z, h, Dh)dx\wedge dy  = B (|z|)\Im\left[\frac{h_T}{h}\right]dx\wedge dy ,\end{equation} where $B \in  L^1(r,R)$.
Thus, for all $h \in \mathcal{H}_j(\mathbb{A}, \mathbb{B})$ we have

\begin{equation}\label{easi1}\int_\A
L(z, h, Dh) dx\wedge dy =\int_r^R
\frac{B(t)}{t}
\int_{|z|=t}
\left(\frac{\partial \mathrm{Arg}\, h}{\partial \tau}
d\tau\right) dt =2\pi j\int_r^R
{B(t)}
dt\end{equation}

The idea of using these free Lagrangians is to establish a general sub-gradient
type inequality for the integrand with two independent parameters: $t \in (r,R)$ and
$s\in (r^\ast,R^\ast)$, for the functions: $g$
as follows:
$$ |g_T|^2+ |g_N|^2\ge  X(s)\frac{|g|_N}{t}
+ Y(t)
\frac{|g|}{s }\mbox{Im}\frac{g_T}{g} + Z(s)\det\left[ Dg\right] + W(t),$$
where the coefficients $X$ and $Z$ are functions in the interval $(r^\ast,R^\ast)$, while $Y$ and
$W$ are functions in the interval $(r,R)$. We will choose $t = |z|$ and $s = |g(z)|$ to
obtain the corresponding free Lagrangians on the right-hand side of the above inequality.

Suppose $g\in H_j(\A, \B)$ and consider the following formula for Dirichlet energy
$$\mathcal{E}[g]=\int_{\A}( |g_N|^2+|g_T|^2) dxdy.$$ Here for $z=\rho e^{\imath\tau}$,  $$g_N(z) =\frac{\partial g(z)}{\partial \rho} $$ and $$g_T =\frac{1}{\rho }\frac{\partial g(z)}{\partial \tau}. $$

Then the Jacobian of a mapping $g$ can be written as $$\det [Dg(z)] = \Im (g_T\overline{g_N}).$$

\section{The proof of Theorem \ref{theoj}}
Since we already proved the existence of radial minimizers, it remains to prove inequality and equality statement. The proof consists of three parts.
\subsection{Under Nitsche bound}

Let us consider three subcases $c_1=0$, i.e. the case $R=r^j$, $c_1>0$ and $c_1<0$.\label{sec-5-2}
In what follows, we let $g(z)=(g_1(z), g_2(z))$ to be the $j-$degree mapping between annuli $\A$ and $\B$ mapping inner/outer boundary onto the inner/outed bounrary.
\begin{description}
  \item[(i)] Conformal  case, $c_1=0$, i.e. $R=r^j$.
\end{description}
From the following trivial inequality $$|g_N|^2+|g_T|^2\ge \det(Dg(z)),$$ it follows that $$\mathcal{E}[g]\ge \int_{\A}\det(Dg(z)) dxdy.$$

Now because $g$ is a $j-$degree mapping we obtain the inequality $$\mathcal{E}[g]\ge j |\B |=\int_{\A}\det(Dg^\circ(z))(z)dxdy.$$

Namely in this case $g^\circ(z) = z^j$ and so $\det(Dg^\circ(z))(z)=j^2 |z|^{2j-2}$. Thus $$\int_{\A}\det(Dg^\circ(z))(z)dxdy=2\pi j^2 \int_1^r \rho^{2j-1}d\rho=\pi j (r^{2j}-1)= j|\B|.$$
\begin{description}
  \item[(ii)] Elastic case, $c_1>0$, i.e. $R>r^j$.
\end{description}

For every $\alpha, \beta\in \mathbb{R}$, the following trivial inequality \begin{equation*}\label{hengy}(\alpha |g_N|-\beta|g_T|)^2\ge 0\end{equation*} is equivalent to
\begin{equation}\label{general} |g_N|^2+|g_T|^2\ge (1- \alpha^2)|g_N|^2 +(1- \beta^2) |g_T|^2 + 2 \alpha \beta   |g_T| |g_N|.\end{equation}
If we take $\beta=1$, then it reduces to the following simple inequality
\begin{equation*}\label{first} |g_N|^2+|g_T|^2\ge (1- \alpha^2)|g_N|^2 +2  \alpha  |g_T| |g_N|,\end{equation*}
where the equality holds if and only if
\begin{equation}\label{simpli11eq1}\alpha|g_N|=|g_T|.\end{equation}
Furthermore, the following inequality holds:
\begin{equation}\label{ess} (1- \alpha^2)|g_N|^2 +2  \alpha  |g_T| |g_N|\ge 2( 1- \alpha^2)\gamma {|g_N|}+ 2 \alpha \det(Dg(z))  -{(1 - \alpha^2) \gamma^2},\end{equation}
where $\gamma$ is an arbitrary constant.
It is easy to see that the equality is attained in \eqref{ess},  if $\gamma=|g_N|$.

Let  $t=|z|$, $s=|g(z)|$. Suppose
\begin{equation*}\label{b}\alpha=\frac{js}{\sqrt{c_1+j^2s^2}},\end{equation*}
\begin{equation*}\label{bast}\alpha_\ast=\frac{jG(t)}{\sqrt{c_1+j^2G^2(t)}},\end{equation*} and
\begin{equation*}\label{bcapital}A=\dot G(t)\sqrt{\frac{ 1-{\alpha^2_\ast}}{1-\alpha^2 }},\end{equation*}
where $G$ is a solution to the differential equation \eqref{pm}
Moreover, in view of \eqref{fprims} \begin{equation}\label{lambda}\begin{split}
|Dg|^2&\ge 2\left[( 1-\alpha^2 )A \right]{|g_N|}+\frac{2 jsF'(s)}{F(s)}\det[Dg] -(1- \alpha^2 )A^2
\\ &\ge 2\left[( 1-\alpha^2 )A \right]{|g|_N}+2j\frac{sF'(s)}{F(s)}\det[Dg] -(1- \alpha^2 )A^2
.\end{split}\end{equation}
Furthermore, in view of \eqref{dgs},
\begin{equation}\label{ss}\begin{split}\gamma(s)&\stackrel{\text{def}}{=} 2t( 1- \alpha^2)A \\&=2t \sqrt{(1- {\alpha^2_\ast})(1-  \alpha^2)} \dot G(t)\\&=\frac{2c_1t\dot G(t)}{\sqrt{j^2G^2(t)+c_1}}\cdot\frac{1}{\sqrt{j^2s^2+c_1}}\\&=2c_1\frac{1}{\sqrt{j^2s^2+c_1}},\end{split}\end{equation} and
\begin{equation}\label{tt}\begin{split}
\delta(t)&\stackrel{\text{def}}{=}-(1- \alpha^2)A^2\\&=-\dot G^2(t)(1- \alpha^2_\ast)\\&=-c_1\frac{\dot G^2(t)}{\sqrt{c_1+j^2G(t)^2}}\\&=-\frac{c_1}{t^2}.\end{split}\end{equation}

So for $t=|z|$ and $s=|g(z)|$, by using the relations \eqref{ess}, \eqref{lambda}, and \eqref{ss}, we obtain \begin{equation*}\begin{split}
\mathcal{E}[g]&=\int_{\A} ( |g_N|^2+|g_T|^2)dx\wedge dy\\&\ge \int_{\A} \left\{2\left[( 1-\alpha^2)A \right]{|g|_N}+2j\frac{sF'(s)}{F(s)}\det\left[ Dg\right] -(1- \alpha^2)A^2\right\} dx\wedge dy
\\&=\int_{\A}  \gamma(s) \frac{|g|_N}{t} dx\wedge dy+2j\int_{\A}  \frac{|g(z)|F'(|g(z)|)}{F(|g(z)|)}\det\left[ Dg\right] dx\wedge dy+\int_{\A}\delta(t)dx\wedge dy.
\end{split}
\end{equation*}
By using the formulas \eqref{easi0}, \eqref{easi}, \eqref{easi1}, and \eqref{ss} we obtain
\begin{equation*}\label{upsi}\mathcal{E}[g]\ge 2\pi\int_1^R \gamma(s) ds+4\pi j^2 \int_1^R   \frac{s^2F'(s)}{F(s)}ds+2\pi \int_{1}^r t\delta(t)dt. \end{equation*}
Then
\begin{equation*}\label{upsi2}\mathcal{E}[g]\ge 4\pi c_1\int_1^R \frac{1}{\sqrt{j^2s^2+c_1}} ds+4\pi j^2 \int_1^R   \frac{s^2F'(s)}{F(s)}ds-2\pi c_1 \int_{1}^r\frac{1}{t}dt. \end{equation*}
According to \begin{equation}\label{okey}\int_{1}^{G(t)} \frac{1}{\sqrt{j^2s^2+c_1}} \, ds=\log t,\end{equation} we see that
\begin{equation*}\label{okey2}\begin{split}\mathcal{E}[g]&\ge 4\pi c_1\log r-2\pi c_1 \log r+4\pi j^2\int_1^R   \frac{s^2F'(s)}{F(s)}ds\\&=2\pi c_1\log r+4\pi j^2\int_1^R   \frac{s^2F'(s)}{F(s)}ds.\end{split} \end{equation*}
Changing the variable $t=F(s)$ in the last integral, we get
\begin{equation}\label{previous}
\int_1^R   \frac{s^2F'(s)}{F(s)}ds=\int_1^r \frac{G^2(t)dt}{t}.
\end{equation}
In view of \eqref{energyof}, we see from \eqref{previous} that $$\mathcal{E}[g]\ge 2\pi c_1\log r+4\pi j^2\int_1^r \frac{G^2(t)dt}{t}=\mathcal{E}[g^\circ].$$
This finishes the proof of this case up to the equality statement.

\begin{description}
  \item[(iii)] Non-elastic case, $-j^2\le c_1< 0$,  i.e. $R<r^j$
\end{description}

Suppose $g$ is a smooth $j-$degree mapping  between two annuli $\A$ and $\B$. Let $z\in \A$, $g(z)\in \B$ and let $t=|z|$, $s=|g(z)|$. If we put $\alpha=1$ in \eqref{general}, then we have the following inequality
\begin{equation}\label{simpli11} |g_N|^2+|g_T|^2\ge (1- \beta^2)|g_T|^2 +2 \beta   |g_T| |g_N|,\end{equation}
where the equality holds if and only if
\begin{equation}\label{simpli11eq}\beta|g_T|=|g_N|.\end{equation}

%Since $g$ maps the outer boundary into the outer boundary, it follows that $\det[Dg]>0$.
Furthermore,
\begin{equation*}\label{equat}\det[Dg]=\mbox{Im}(g_T\overline{g_N})=|g_T \overline{g_N}|\end{equation*}
if and only if
\begin{equation}\label{imaine}g_T\overline{g_N}\in \imath\mathbb{R}.\end{equation}
We continue with the following inequality
\begin{equation*}\label{simpli1} |g_N|^2+|g_T|^2\ge \mathcal{R}:{=} 2( 1- \beta^2)B {|g_T|}+ 2 \beta \det(D g)-{( 1- \beta^2) B^2},\end{equation*}
where $B$ is a constant.

Since $\mbox{Im}\left[\frac{g_T}{g}\right]\le\left[\frac{|g_T|}{|g|}\right]$, we infer that
\begin{equation*}\label{imi}|g_T|\ge |g|\mbox{Im}\left[\frac{g_T}{g}\right].\end{equation*}
The above equality is attained if and only if
\begin{equation}\label{pavaru}
g_T\overline{g}\in \imath\mathbb{R}.
\end{equation}
Now
\begin{equation*}\label{help}
\mathcal{R}\ge [2( 1- \beta^2)B]{|g|\mbox{Im}\left[\frac{g_T}{g}\right]}+ 2\beta\det\left[ Dg\right] -{(1 - \beta^2) B^2}.
\end{equation*}
The above equality is attained if and only if the conditions \eqref{imaine}, \eqref{simpli11eq}, \eqref{pavaru} and the following equality
\begin{equation}\label{pavar} B=|g_T| \end{equation}
are satisfied.

Moreover, these conditions imply that $g$ is radial. Now, we choose
\begin{equation}\label{beta}\beta=\frac{\sqrt{j^2s^2+c_1}}{js}\end{equation}
and
\begin{equation}\label{AA}B=\frac{j s}{t}.\end{equation}
Then, for $t=|z|$ and $s=|g(z)|$,  we get
\begin{equation}\label{pr}\begin{split}-{(1 - \beta^2) B^2}&=\mu(t)\stackrel{\text{def}}{=} \frac{c_1}{t^2}.\end{split}\end{equation}
and
\begin{equation}\label{pri}\begin{split}2s( 1- \beta^2)B=\nu(t)\stackrel{\text{def}}{=} -\frac{2c_1}{jt}.\end{split}\end{equation}
By using the relations \eqref{simpli11}, \eqref{pr}, and \eqref{pri} we obtain \begin{equation}\label{922}\begin{split}
\mathcal{E}[g]&=\int_{\A} ( |g_N|^2+|g_T|^2)dx\wedge dy
\\&\ge \int_{\A} \bigg[[2(1 - \beta^2)B ]{{|g|\mbox{Im}\left[\frac{g_T}{g}\right]}}\\&+ 2 \frac{\sqrt{j^2s^2+c_1}}{js} \det\left[ Dg\right] -{(1 - \beta^2) B^2}\bigg] dx\wedge dy
\\&=\int_{\A}  \nu(t) \Im\left[\frac{g_T}{g}\right] dx\wedge dy\\&+2\int_{\A}  \frac{\sqrt{j^2s^2+c_1}}{js} \det\left[ Dg\right] dx\wedge dy+\int_{\A}\mu(t)dx\wedge dy.
\end{split}
\end{equation}

By using the formulas \eqref{easi0}, \eqref{easi}, \eqref{easi1}, and \eqref{922} we obtain

\begin{equation*}\label{upsi21}\mathcal{E}[g]\ge 2\pi j\int_1^r  \nu(t) dt+4\pi \int_1^R \sqrt{j^2 s^2+c_1 }ds+2\pi \int_{1}^r t\mu(t)dt. \end{equation*}
Now we use the equality
$$\sqrt{j^2s^2+c_1 }=\frac{j^2s^2+c_1}{\sqrt{j^2s^2+c_1}}$$ to deduce that
$$\int_{1}^{R}\sqrt{j^2s^2+c_1 }ds=I_1+I_2=c_1\int _{1}^{R}\frac{1}{\sqrt{j^2s^2+c_1}}ds+\int_{1}^{R}j^2s^2\frac{1}{\sqrt{j^2s^2+c_1}}ds.$$
By using \eqref{okey} we get $$I_1=c_1\log r.$$
By using  change of variables $t=F(s)$, in view of
\begin{equation}\label{kina}F(\tau)=r\exp\left[\int_{1}^{\tau} \frac{1}{\sqrt{s^2+c_1}}ds\right],\end{equation}  we get $$I_2=j^2\int_{1}^{R}s^2\frac{F'(s)}{F(s)}ds=j^2\int_1^r\frac{G^2(t)}{t}dt.$$
Moreover, we have $$2\pi j\int_1^r  \nu(t) dt=-4\pi c_1\log r.$$
We also have $$2\pi \int_{1}^r t\mu(t)dt= 2\pi c_1 \log r.$$
Summing all those formulas and again from \eqref{energyof}, we see that $$\mathcal{E}[g]\ge  2\pi c_1\log r+4\pi j^2\int_1^r \frac{G^2(t)dt}{t}=\mathcal{E}[g^\circ].$$

\subsection{Below Nitsche bound}

In this case we assume that $\A=\A(\rho,1]\cup \A[1,r)$. We repeat the previous argument (the case $\mathrm{(iii)}$). In this case $c_1=-j^2$ and  we have \begin{equation}\label{upsi212}\mathcal{E}[g]\ge 2\pi j\int_\rho^r  \nu(t) dt+4\pi \int_1^R \sqrt{j^2 s^2-j^2 }ds+2\pi \int_{\rho}^r t\mu(t)dt. \end{equation}
Further as before we obtain $$4\pi \int_{1}^{R}\sqrt{j^2s^2-j^2 }ds= 4\pi j^2\int_1 ^r \frac{G^2(t)dt}{t}-4\pi j^2 \log r,$$
and in view of \eqref{pr} and \eqref{pri}
\[\begin{split}\mathcal{E}[g]&\ge 4\pi j^2\log \frac{r}{\rho}+4\pi j^2\int_1 ^r \frac{G^2(t)dt}{t}-4\pi j^2 \log r-2\pi j^2 \log \frac{\rho}{r}
\\&=4\pi j^2 \int_1^r\frac{G^2(t)}{t}dt-2\pi j^2\log r-2\pi j^2\log \rho\\&=\mathcal{E}[g^\diamond].\end{split}\]
\subsection{The equality part}
To finish the proof of the theorem, we need to consider the equality case. Consider only the non-elastic case, and notice that the proof of the remaining case is similar. Assume that the equality is attained in our theorem for a certain mapping $g$. Let $z=t e^{\imath\tau}$, and write $g(t,\tau)=g(z)$. Assume that $g(t,\tau)= u(t,\tau) e^{\imath v(t,\tau)}$ for some real functions $u$ and $v$, with $u(1,\tau)=1$, $u(r,\tau)=R$, $u(t,\tau+2\pi)=u(t,\tau)$, $v(t,\tau+2\pi)=v(t,\tau)$ for every $t,\tau\in [1,r]\times \mathbb{R}$. Then from \eqref{imaine} and \eqref{pavaru}, we see that $$\frac{\partial_t g(t,\tau)}{g(t,\tau)}=\imath \partial_t v(t,\tau)+\frac{\partial_t u(t,\tau)}{u(t,\tau)}\in \mathbb{R}.$$ Since $g$ is $j-$degree mapping it follows that  $v(t,\tau)=\phi(\tau)$ for a certain $j-$degree mapping $\phi$ of the unit circle onto itself. Now from \eqref{imaine}, we have $$\frac{\left(\imath u(t,\tau) \phi'(\tau)+\partial_\tau u(t,\tau)\right)}{t}\in \imath \mathbb{R}.$$ Thus $u(t,\tau)= G(t)$. Now from \eqref{simpli11eq} we see that
$$\beta \frac{G(t) \phi'(\tau)}{r}=G'(t),$$ where, as in \eqref{beta}, \begin{equation}\label{beta1}\beta = \frac{\sqrt{j^2s^2+c_1}}{js}=\frac{\sqrt{G(t)^2+c_1}}{j G(t)}.\end{equation}

This implies that $\phi'(\tau)=\mathrm{const}=j$, and thus $\phi(\tau) =j \tau +c$. Now from \eqref{simpli11eq} and \eqref{beta1} we see that $G$ satisfies the equation \eqref{pm}, which means that $g$ is a $j-$degree harmonic mapping.

%\section{Another paper-Existence of $j-$degree mappings between certain double-connected domains??? }

%{\bibliographystyle{abbrv} \bibliography{references}}

\end{document}